\documentclass[submission,copyright]{eptcs}
\providecommand{\event}{GandALF 2019}
\usepackage{breakurl}            
\usepackage{underscore}          
\usepackage[T1]{fontenc} 
\usepackage{amsmath}
\usepackage{amssymb}
\usepackage{amsthm}
\usepackage{hyperref}
\usepackage{float}
\usepackage{tikz}
\usetikzlibrary{arrows,shapes}

\newcommand{\eps}{\varepsilon}

\newcommand{\sgn}{\mathop{\mathrm{sgn}}}

\newcommand{\calG}{\ensuremath{\mathcal{G}}}

\newcommand{\Reach}{\ensuremath{\mathsf{Reach}}}
\newcommand{\Safe}{\ensuremath{\mathsf{Safe}}}

\title{A Stay-in-a-Set Game without a\\ Stationary Equilibrium}
\def\titlerunning{A Stay-in-a-Set Game without a Stationary Equilibrium}
\author{
Kristoffer Arnsfelt Hansen
\institute{Aarhus University\\Department of Computer Science}
\email{arnsfelt@cs.au.dk}
\and
Mikhail Raskin
\thanks{
The author has received funding from the European Research Council (ERC) under the European Union’s Horizon 2020 research and innovation programme under grant agreement No 787367 (PaVeS) and from the French National Research Agency (ANR project GraphEn / ANR-15-CE40-0009).
}
\institute{Technical University of Munich\\Department of Informatics}
\email{raskin@mccme.ru}
}
\def\authorrunning{K.~A.~Hansen \& M.~Raskin}
\begin{document}
\maketitle

\begin{abstract}
  We give an example of a finite-state two-player turn-based
  stochastic game with safety objectives for both players which has no
  stationary Nash equilibrium. This answers an open question of Secchi
  and Sudderth.
\end{abstract}

\section{Introduction}
Stochastic games provide a general model for studying dynamic
interactions between players whose actions affect the state of the
environment. The change in state is described by a probability
distribution called the law of motion. The first such games were
introduced by Shapley~\cite{PNAS:Shapley53}. We may view his model as
discrete-time finite two-player zero-sum games, where players receive
immediate payoff in each round of play and \emph{discount} future
payoffs. Shapley proved that in such games optimal stationary (or,
memoryless) strategies exists. The initial model of Shapley has since
been extensively extended and studied in many variations. With each
model the main question is existence of optimal strategies or a Nash
equilibrium. Next, it is of interest how complicated such strategies
must be. We shall limit our discussion to discrete-time games having
an arbitrary but finite number of states.

Everett~\cite{AMS:Everett57} defined \emph{recursive games}, where
players only receive a (possibly) non-zero payoff when play terminates
by entering special absorbing states. These payoffs are also called
terminal payoffs. While players are no longer guaranteed to have
optimal strategies, Everett proved that they do have $\eps$-optimal
stationary strategies. Gillette~\cite{AMS:Gillette1957} considered
finite two-player zero-sum games where players again receive immediate
payoffs each round of play, but now evaluate their payoff as the
average of the immediate payoffs received (limit average payoff). Here
players are no longer guaranteed to have $\eps$-optimal stationary
strategies, but as shown by Mertens and
Neyman~\cite{IJGT:MertensNeyman1981} they do have $\eps$-optimal
strategies. An even more general result was obtained by
Martin~\cite{JSYML:Martin98} showing that for two-player zero-sum
games where payoffs are Borel measurable functions of the history of
play, the players have $\eps$-optimal strategies. Here the extension
from deterministic games (i.e.\ games having a deterministic law of
motion) to the general case is due to an observation of Maitra and
Sudderth~\cite{IJGT:MaitraS98}.

For non-zero sum games much less is known. For discounted payoffs, a
Nash equilibrium exists in stationary strategies as shown by
Fink~\cite{Fink1964} and Takahashi~\cite{Takahashi1964}. The existence
of $\eps$-Nash equilibrium in recursive games is an open problem, even
for~three players. In addition, Flesch, Thuijsman and
Vrieze~\cite{MOR:FleschTV96} gave an example of a two-player recursive
game without \emph{stationary} $\eps$-Nash equilibrium.
Vieille~\cite{IJM:Vieille2000a,IJM:Vieille2000b} proved existence of
$\eps$-Nash equilibrium in every two-player game with limit-average
payoff.

Mertens and Neyman (cf.\ \cite{Mertens87}) showed, using the
celebrated determinacy result by Martin~\cite{AM:Martin75}, that an
$\eps$-Nash equilibrium exists in any \emph{turn-based} (i.e.\ perfect
information) game with Borel payoff functions. Later this was observed
again by Chatterjee~et~al.~\cite{CSL:ChatterjeeMJ04}. When the payoff
function has finite range, an actual Nash equilibrium exists. This is
particularly the case of deterministic games where the payoff function
is the indicator function of a Borel set. We refer to the indicator
function of a Borel set as well as the set itself as a Borel
\emph{winning set} or Borel \emph{objective}.

The most basic of these are given by the open and closed sets. Given a
set of states $T$, the \emph{reachability objective} $\Reach(T)$ given
by $T$ consists of the histories of play that visit a state in
$T$. The \emph{safety objective} $\Safe(T)$ given by $T$ consists of
the histories of play that stays within the states in $T$. These
winning sets are the open and closed Borel objectives typically
studied, and they have applications in the verification and synthesis
of reactive systems~\cite{JCSS:ChatterjeeH12}.

Games where the players have reachability or safety objectives are
closely related to recursive games. First note that for a given
recursive game, after normalizing all payoffs to be in the range
$[-1,1]$, every terminal payoff vector can be written as a convex
combination of payoff vectors having only entries from the set
$\{-1,0,1\}$. This means that any absorbing state can be replaced by a
set of absorbing states where all players have payoffs in the set
$\{-1,0,1\}$ as well by modifying the (probabilistic) law of motion
accordingly. Then, if a player only receive terminal payoffs from the
set $\{-1,0\}$, this is equivalent to a safety objective, and likewise
if a player only receive terminal payoffs from the set $\{0,1\}$, this
is equivalent to a reachability objective.

Secchi and Sudderth~\cite{IJGT:SecchiS02} considered the class of
games where each player has a safety objective, and called these games
for \emph{stay-in-a-set games}. For these games they proved existence
of a Nash equilibrium in any (finite) stay-in-a-set game. The
equilibrium strategies are not stationary but prescribe, as a function
of the set of the players whose safety objective has not yet been
violated, a stationary strategy profile. 
This means that for $n$ players just $n$ bits of (shared) memory
are needed for implementing a Nash equilibrium.
A natural open question
raised by Secchi and Sudderth was then existence of a stationary Nash
equilibrium. We give an example of a two-player game without a
stationary Nash equilibrium. Our game is furthermore turn-based. By
example we also illustrate the Nash equilibria obtained from the proof
of Secchi and Sudderth. They rely crucially on the willingness of the
second player to change strategy after already having lost.
In our example it is not necessary to remember whether the first 
player has lost, so a single bit of shared memory is enough for
implementing the equilibrium.
Finally we
note that players do have a stationary $\eps$-Nash equilibrium.

It is necessary that our example game is not deterministic.  In fact,
in every deterministic two-player turn-based games, where each player
has a reachability or a safety objective, a Nash equilibrium exists in
positional (i.e.\ pure and memoryless) strategies. This follows from
the fact that two-player zero-sum games with a reachability and safety
objective are positionally determined. Thus in the non-zero sum game
it is either the case that one of the two players may guarantee a win
(and relative to that we let the other player play optimally) or it is
the case that both players can ensure that the opponent loses.
\enlargethispage{\baselineskip}

\section{The game}
The game $\calG$ we consider is played by two players each taking
turns in choosing whether to continue the game or to attempt to quit
the game, with Player~1 making the first choice. A choice of the quit
action by one of the players is successful with probability~$3/4$, and
otherwise the game continues with the other player as before.  If
Player~1 makes the choice to quit both players win with
probability~$1/2$ and both players lose with the remaining
probability~$1/4$. If Player~2 makes the choice to quit both players
win with probability~$1/4$ and both players lose with the remaining
probability~$1/2$. Finally, Player~2 is incentivized to choose quit by
having the continue action of Player~2 lead to a loss for Player~2
with probability~$1/4$. Infinite play leads to Player~1 winning (and
Player~2 losing with probability~1). This leads to a discontinuity in
the payoff function of Player~1, which is crucial for our example.

The game $\calG$ is illustrated in Figure~\ref{FIG:Game} and is
modeled with a set of~5 states $\{1,2,W,L,L2\}$, with Player~1
controlling state~1, Player~2 controlling state~2. State~$L2$ exists
merely to enforce a loss to Player~2, whereas the states $W$ and $L$
are winning and losing states of both players, respectively. The game
$\calG$ is a stay-in-a-set game with the safe sets of the two players
being $G_1=\{1,2,W,L2\}$ and $G_2=\{1,2,W\}$, respectively.  The
diamond-shaped nodes in Figure~\ref{FIG:Game} are used to indicate the
probabilistic transitions.

A stationary strategy profile $\sigma$ in $\calG$ can be described by
a pair of probabilities $(p_1,p_2)$, where $p_i$ is the probability that
Player~$i$ chooses the quit action $q$, when in state~$i$ (and thus
$(1-p_i)$ is the probability that Player~$i$ chooses the continue
action $c$, when in state~$i$).

\begin{figure}[H]
 \centering
\begin{tikzpicture}[>=stealth',->,bend angle=45,auto,node distance=2cm]
\tikzstyle{p1node}=[regular polygon,regular polygon sides=4,thick,draw=black,minimum size=8mm]
\tikzstyle{p2node}=[regular polygon,regular polygon sides=8,thick,draw=black,minimum size=8mm]
\tikzstyle{othernode}=[circle,thick,draw=black,minimum size=8mm]
\tikzstyle{randomnode}=[diamond,thick,draw=black]
\node [p1node] (p1) at (0,0) {$1$};
\node [p2node] (p2) at (8,0) {$2$};
\node [othernode] (L2) at (4,0) {$L2$};
\node [randomnode] (p1quit) [below of=p1] {};
\node [randomnode] (p2cont) [left of=p2] {};
\node [randomnode] (p2quit) [below of=p2] {};
\node [othernode] (W) at (4,-2) {$W$};
\node [othernode] (L) at (4,-3) {$L$};

\draw[thick] [rounded corners] (p1)|- node[near start]{$c$} (2.5,2) -|(p2);
\draw[thick] (p1) -- node[swap]{$q$} (p1quit);
\draw (p1quit) -- node[pos=0.2]{$1/4$} (p2);
\draw (p1quit) -- node[pos=.6]{$1/2$}(W);
\draw (p1quit) -- node[swap,pos=.4]{$1/4$}(L);

\draw[thick] (p2) -- node[swap]{$c$} (p2cont);
\draw[thick] (p2) -- node{$q$} (p2quit);
\draw[thick] (L2) -- (p1);
\draw (p2cont) -- node[swap]{$1/4$} (L2);
\draw [rounded corners] (p2cont) |- node[swap,pos=0.615]{$3/4$} (2,1) -- (p1);
\draw (p2quit) -- node[swap,pos=0.2]{$1/4$} (p1);
\draw (p2quit) -- node[swap,pos=.6]{$1/4$}(W);
\draw (p2quit) -- node[pos=.4]{$1/2$}(L);

\draw[thick] (W) edge [loop above] (W);
\draw[thick] (L) edge [loop below] (L);

\draw (-1,0) -- (p1);
\end{tikzpicture}
\caption{The game $\calG$.}
\label{FIG:Game}
\end{figure}

\subsection{No stationary Nash equilibrium}
We give here a simple analysis showing that no stationary Nash
equilibrium exists in $\calG$. We can place all plays of $\calG$ in~3
groups. Group~1 are plays where Player~1 quits successfully, group~2
are plays where Player~2 quits successfully, and group~3 are plays
that never reach $W$ or $L$.

Consider a stationary strategy profile $\sigma$ given by $(p_1,p_2)$.
When $p_1=p_2=0$ the play belongs to group~3, where Player~1 wins and
Player~2 loses with probability~$1$. When $p_1>0$ or $p_2>0$ the
play belongs to group~1 or group~2 with probability~$1$. The players
both prefer a play from group~1 where Player~1 is the player to quit
successfully.

Suppose that $\sigma$ is a Nash equilibrium. If $p_2=0$, then the only
best reply of Player~1 is to have $p_1=0$, since otherwise $L$ is
reached with positive probability. But if also $p_1=0$, Player~2 loses
with probability~1, whereas $p_2>0$ would lead to reaching $W$ with
positive probability. This rules out having $p_2=0$ in a Nash
equilibrium.

Suppose now $p_2>0$, which means that the play belongs to group~$1$ or
group~$2$ with probability~1. The probability that the play belongs to
group~1 strictly increases with $p_1$, and it follows that we must
have $p_1=1$. But this is also not a Nash equilibrium, as Player~2
would then be better off having $p_2=0$. Indeed, let us consider a
play from state~2 until the play either returns to state~2, reaches
state~$W$ before returning to state~2, or reaches state~$L$ or
state~$L2$ before returning to state~2. We denote these events a
return, a win, or a loss.

The quit action for Player~2 has probability
$1/2 + (1/4)(1/4) = 9/16$ of a loss, probability
$1/4+(1/4)(1/2) = 6/16$ of a win and $(1/4)(1/4)=1/16$ of
a return. The continue action has probability
$1/4+(3/4)(1/4)=7/16$ of a loss, probability
$(3/4)(1/2)=6/16$ of a win and $(3/4)(1/4)=3/16$ of
a return. Since a return is better than a loss for Player~2, this rules out
$p_2>0$ in a Nash equilibrium as well.

\subsection{Detailed payoff analysis}

For $i,j \in \{1,2\}$, let $u_{i,j}=u_{i,j}(p_1,p_2)$ be the payoff to
Player~$i$ of the strategy profile $(p_1,p_2)$ when starting play in
state $j$. The payoffs satisfy the following equations
\begin{align*}
  u_{1,1} & = ((1-p_1)+p_1/4)u_{1,2}+p_1/2  = (1-3p_1/4)u_{1,2}+p_1/2\\
  u_{1,2} & = ((1-p_2)+p_2/4)u_{1,1}+p_2/4  = (1-3p_2/4)u_{1,1}+p_2/4\\
  u_{2,1} & = ((1-p_1)+p_1/4)u_{2,2}+p_1/2  = (1-3p_1/4)u_{2,2}+p_1/2\\
  u_{2,2} & = (3/4(1-p_2)+p_2/4)u_{2,1}+p_2/4  = (3/4-p_2/2)u_{2,1}+p_2/4\\
\end{align*}
and from these follows further
\begin{align*}
  u_{1,1} & = (1-3p_1/4)(1-3p_2/4)u_{1,1} + (1-3p_1/4)p_2/4 + p_1/2\\
  u_{2,1} & = (1-3p_1/4)(3/4-p_2/2)u_{2,1}+ (1-3p_1/4)p_2/4 + p_1/2
\end{align*}
When both $p_1=0$ and $p_2=0$ we have that $u_{1,1}=u_{1,2}=1$ and
$u_{2,1}=u_{2,2}=0$.  When at least one of $p_1>0$ or $p_2>0$ holds,
we can solve for $u_{1,1}$, and likewise we can always solve for
$u_{2,1}$ to obtain
\begin{align*}
  u_{1,1} &=\frac{(1-3p_1/4)p_2/4 + p_1/2}{1-(1-3p_1/4)(1-3p_2/4)}  = \frac{(8-3p_2)p_1+4p_2}{(12-9p_2)p_1+12p_2} \\ &= \frac{(4-3p_1)p_2+8p_1}{(12-9p_1)p_2+12p_1} \\
  u_{2,1} &=\frac{(1-3p_1/4)p_2/4 + p_1/2}{1-(1-3p_1/4)(3/4-p_2/2)} = \frac{(4-3p_1)p_2+8p_1}{(8-6p_1)p_2+9p_1+4} \\ &= \frac{(8-3p_2)p_1+4p_2}{(9-6p_2)p_1+8p_2+4}
\end{align*}
Using
$\frac{\partial}{\partial x} \frac{ax+b}{cx+d}=
\frac{ad-bc}{(cx+d)^2}$, we find that
\[
  \sgn\left(\frac{\partial}{\partial p_1}u_{1,1}\right)=\sgn(48p_2)=1 \enspace \text{ for all } p_2>0 \enspace ,\\
\]
And likewise $\sgn\left(\frac{\partial}{\partial  p_2}u_{2,1}\right)  =\sgn((4-3p_1)(4-7p_1))$, which means that
\[
\sgn\left(\frac{\partial}{\partial  p_2}u_{2,1}\right) = \begin{cases}\phantom{-}1 & \text{if } p_1<\frac{4}{7}\\\phantom{-}0 & \text{if } p_1=\frac{4}{7}\\-1 & \text{if } p_1>\frac{4}{7} \end{cases} \enspace .
\]

The function $u_{2,1}$ is continuous in the entire domain, whereas the
function $u_{1,1}$ has a single discontinuity when $p_1=p_2=0$. Note
that $u_{1,1}(p_1,0)=\frac{2}{3}$ for all $p_1>0$, and
$u_{1,1}(0,p_2)=\frac{1}{3}$ for all $p_2>0$.

The best replies of the players are as follows. If $p_2=0$, the only
best reply of Player~1 is to have $p_1=0$, giving $u_{1,1}=1$. If
$p_2>0$, the only best reply of Player~1 is to have $p_1=1$, giving
$u_{1,1}=\frac{8+p_2}{12+3p_2}$. If $p_1<4/7$, the only best reply for
Player~2 is to have $p_2=1$, giving
$u_{2,1}=\frac{4+5p_1}{12+3p_1}$. When $p_1=4/7$, Player~2 has no
preferred action. Finally, if $p_1>4/7$ the only best reply of
Player~2 is to have $p_2=0$, giving $u_{2,1}=\frac{8p_1}{4+9p_1}$.

\subsection{Nash equilibria}
We give here two examples of Nash equilibria in the game following the general result of
Secchi and Sudderth~\cite{IJGT:SecchiS02}. The idea is that once
Player~2 has lost by entering state $L2$ the incentive of
Player~2 is removed and all strategies are equally good.

Suppose first that Player~2 commits to always playing the continue
action after entering state $L2$. The best reply of Player~1 is then
to always play the continue action as well, ending up with
payoff~1. We may thus consider the modified game $\calG'$ that stops
when entering $L2$ upon which Player~1 receives payoff~1. This lead to
the modified equation
\[
  u_{1,2} = (3/4(1-p_2)+p_2/4)u_{1,1}+1/4  = (3/4-p_2/2)u_{1,1}+1/4
\]
giving
\[
  u_{1,1} = (1-3p_1/4)(3/4-p_2/2)u_{1,1}+(1-3p_1/4)/4+p_1/2 
\]
which solves to
\[
  u_{1,1} = \frac{(1-3p_1/4)/4+p_1/2}{1-(1-3p_1/4)(3/4-p_2/2)} = \frac{5p_1+4}{(9-6p_2)p_1+4+8p_2}
\]
and we see that $\sgn\left(\frac{\partial}{\partial p_1}u_{1,1}\right)=\sgn(64p_2-16)$.

A Nash equilibria is thus that the players play the quit action with
probabilities $p_1=4/7$ and $p_2=1/4$ respectively
until state $L2$ is reached and after which both players play the quit
action with probability $p_1=p_2=0$. The equilibrium payoffs are 
$u_{1,1}=2/3$ and $u_{2,1}=1/2$.

Suppose next that Player~2 commits to always playing the quit action
after entering state $L2$. The best reply of Player~1 is then to
always play the quit action as well, ending up with
payoff~$3/5$. The modified game $\calG'$ now has the
equation
\[
  \begin{split}
u_{1,2} = (3/4(1-p_2)+p_2/4)u_{1,1}+3/5(1-p_2)/4+p_2/4 \\ = (3/4-p_2/2)u_{1,1}+3/20+p_2/10
\end{split}
\]
giving
\[
  u_{1,1} = (1-3p_1/4)(3/4-p_2/2)u_{1,1}+(1-3p_1/4)(3/20+p_2/10)+p_1/2 
\]
which solves to
\[
  u_{1,1} = \frac{(1-3p_1/4)(3/20+p_2/10)+p_1/2}{1-(1-3p_1/4)(3/4-p_2/2)} = \frac{(31-6p_2)p_1+8p_2+12}{5(9-6p_2)p_1+8+4p_2}
\]
We find that
$\sgn\left(\frac{\partial}{\partial
    p_1}u_{1,1}\right)=\sgn(1120p_2+80)=1$, which means that the best
reply of Player~1 is always to play the quit action, and in turn the
best reply of Player~2 to that is to always play the continue action.

A Nash equilibria is thus that the players play the quit action with
probabilities $p_1=1$ and $p_2=0$ respectively until state $L2$ is
reached and after which Player~2 changes to playing the quit action
with probability $p_2=1$ as well. The equilibrium payoffs are here
$u_{1,1}=\frac{43}{65}=\frac{2}{3}-\frac{1}{195}$ and
$u_{2,1}=\frac{8}{13} = \frac{1}{2}+\frac{3}{26}$.

\subsection{Stationary $\eps$-Nash equilibrium}
Whereas we have shown that the game $\calG$ has no stationary Nash
equilibrium, it does have $\eps$-Nash equilibria, for any $\eps>0$.

When $\eps<1/3$ no $\eps$-Nash equilibrium can have $p_2=0$. Indeed,
then the only $\eps$-best reply of Player~1 would be the actual best
reply having $p_1=0$. To that, any $\eps$-best reply of Player~2 must
have $p_2>0$, when $\eps<1/3$.

A few examples of $\eps$-Nash equilibria are $\sigma_1$ given by
$p_1=1$ and $p_2=\eps$, $\sigma_2$ given by $p_1=\frac{4}{7}-\eps$ and
$p_2=\eps$, and $\sigma_3$ given by $p_1=\frac{4}{7}+\eps$ and
$p_2=\eps$. We omit the simple task of verifying that these are indeed
$\eps$-Nash equilibria.  In $\sigma_1$ the payoffs are
$u_{1,1}=\frac{2}{3}-O(\eps)$ and $u_{2,1} = \frac{8}{13}-O(\eps)$,
and in both $\sigma_2$ and $\sigma_3$ the payoffs satisfy
$u_{1,1}=\frac{2}{3}-O(\eps)$ and $u_{2,1}=\frac{1}{2}-O(\eps)$.  We
note that Player~1 is playing the best reply in $\sigma_1$,
but is $(\frac{3}{7}-\eps)$-far from the best reply $p_1=1$ in
$\sigma_2$ and $\sigma_3$. Player~2 is playing $\eps$-close to the
best reply $p_2=0$ in $\sigma_1$ and $\sigma_3$, but $(1-\eps)$-far from
the best reply $p_2=1$ in $\sigma_2$.

\section{Conclusion and Further Problems}
We have given a simple example of a two-player turn-based game with
safety objectives for both players without a stationary Nash
equilibrium. A remaining open question is the existence of a
stationary $\eps$-Nash equilibrium when players have safety
objectives, even in the case of two-player turn-based games.

Several related open questions concern games with reachability
objectives or with combinations of reachability and safety
objectives. We first consider the setting where all players have
reachability objectives, also called \emph{reach-a-set
  games}~\cite{CSL:ChatterjeeMJ04}. Flesch, Thuijsman and
Vrieze~\cite{MOR:FleschTV96} give an example of a three-player
recursive game with non-negative payoffs with no stationary
$\eps$-Nash equilibrium. The game is furthermore
deterministic. Simon~\cite{IJM:Simon2006} gave an example of a
two-player recursive game with non-negative payoffs with no stationary
$\eps$-Nash equilibrium. These both give examples of reach-a-set games
without stationary $\eps$-Nash equilibria by the general method of
simulating terminal payoffs with the probabilistic law of motion. The
example of Flesch, Thuijsman and Vrieze is however such that the
terminal payoff vectors satisfy that either none or precisely two
players receive a strictly positive payoff. The payoff vectors where
two players receive strictly positive payoff can (after scaling) be
constructed as unique equilibrium payoffs of win-lose bimatrix
games~\footnote{These payoff vectors are $(1,1)$ and $(3,1)$. It is
  easy to construct two $4 \times 4$ bimatrix games with only payoffs
  from the set $\{0,1\}$ in which the unique equlibrium payofff
  vectors are $(\tfrac{1}{4},\tfrac{1}{4})$ and
  $(\tfrac{3}{4},\tfrac{1}{4})$, respectively, which may replace $(1,1)$ and $(3,1)$}. This then results in
a three-player deterministic reach-a-set game with no stationary
$\eps$-Nash equilibrium.

For two-player games, it was erroneously claimed
(cf.~\cite{FSTTCS:BouyerMS14}) first by
Chatterjee~et~al.~\cite{CSL:ChatterjeeMJ04} and later again by Ummels
and Wojtczak~\cite{CONCUR:UmmelsW11} that a simple adaptation of an
example of a zero-sum game of Everett resulted in a deterministic
reach-a-set game without a Nash equilibrium. Thus it remains an open
question whether every deterministic two-player reach-a-set game has a
Nash equilibrium. It is also an open problem whether every
deterministic two-player reach-a-set game has a stationary $\eps$-Nash
equilibrium. Boros and Gurvich~\cite{MSS:BorosG03} and
Kuipers~et~al.~\cite{EJOR:KuipersFSV09} give an example of a
three-player turn-based recursive game with non-negative payoffs that
has no stationary Nash equilibrium. Do every two-player turn-based
reach-a-set game have a stationary Nash equilibrium?

Little is known when some players have a reachability objective and
some players a safety objectives. In the two-player zero-sum case an
example of Everett~\cite{AMS:Everett57} shows that optimal strategies,
and hence a Nash equilibrium, may fail to exist. On the other hand an
$\eps$-optimal stationary equilibrium always exists. Do every
two-player game where one player has a reachability objective and one
player a safety objetive always have a stationary $\eps$-Nash
equilibrium? In the case of turn-based games, it is an open problem
whether every three-player deterministic game has a stationary Nash
equilibrium. An example given by Boros~et~al.~\cite{DAM:BorosGMOV18}
appears to be close to answer this question. Namely, Boros~et~al.\
construct a three-player deterministic recursive game without a
stationary Nash equilibrium, that may be realized with payoffs such
that player two has only non-negative terminal payoffs and player one
and player three have only non-positive terminal payoffs.

\bibliographystyle{eptcs} \bibliography{SafetyNoNash}

\end{document}